\documentclass[a4paper,12pt,leqno]{article} 

\usepackage{lineno}

\usepackage{amsmath,amsfonts} 
\usepackage{graphics}
\usepackage[T1]{fontenc}
\usepackage{lineno}
\usepackage[utf8]{inputenc}
\usepackage{amsmath,mathrsfs,amssymb,amsfonts,nicefrac}
\usepackage[english]{babel}
\usepackage[active]{srcltx}
\usepackage{tikz}
\usepackage{dsfont}
\usetikzlibrary{patterns}

\numberwithin{equation}{section}

\newtheorem{thm}{Theorem}[section]   
   
\newtheorem{prop}[thm]{Proposition}   
\newtheorem{lm}[thm]{Lemma}

\newcommand{\RR}{\mathbb{R}}   
\newcommand{\R}{\mathbb{R}}   
\newcommand{\di}{\displaystyle}

\newcommand{\NN}{\mathbb{N}}

\renewcommand{\epsilon}{\varepsilon}
\newcommand{\e}{\varepsilon}
\def\un{{\mathbf{1}}}

\begin{document}   
   
\title{\textbf{The shape of a free boundary driven by a line of fast diffusion}}

\author{
{\bf Luis A. Caffarelli} \\ 
The University of Texas at Austin\\
Mathematics Department RLM 8.100, 2515 Speedway Stop C1200\\
 Austin, Texas 78712-1202, U.S.A.  \\
\texttt{caffarel@math.utexas.edu}
\\[2mm]
{\bf Jean-Michel Roquejoffre} \\ 
Institut de Math\'ematiques de Toulouse (UMR CNRS 5219) \\ 
Universit\'e Toulouse III,
118 route de Narbonne\\
31062 Toulouse cedex, France \\ 
\texttt{jean-michel.roquejoffre@math.univ-toulouse.fr}}
  \date{}
\maketitle  

\smallskip 
\begin{center}
{\it  We would like to dedicate this article to Sandro Salsa.  He is a wonderful person and an exceptional, generous mathematician.
We have greatly enjoyed his work and his friendship}
\end{center}

\smallskip
\begin{abstract} \noindent We complete the description, initiated in \cite{CafR}, of  a free boundary travelling at constant speed in a half plane, where  the propagation is controlled  by a line having a large diffusion on its own.
The main result of this work is that the free boundary is asymptotic to a line at infinity, whose angle to the horizontal is dicatated by the velocity of the wave imposed by the line. This helps understanding some rather counter-intuitive numerical simulations of \cite{ACC}.
\end{abstract}

\section{Introduction} 
\label{s1}   
\subsection{Model and question}
The system under study involves  an unknown real $c>0$, an unknown function $u(x,y)$ defined in $\RR^2_-:=\RR_+\times\RR_-$, and an unknown  curve $\Gamma\subset\RR^2_-$ satisfying
 \begin{equation}
\label{e1.2}
\left\{
\begin{array}{rll}
-d\Delta u+c\partial_x u=&0\quad(x,y)\in\{u>0\}\\
\vert\nabla u\vert=&1\quad((x,y)\in\Gamma:=\partial\{u>0\}\\
\ \\
-Du_{xx}+c\partial_xu+1/\mu u_y=&0\quad\hbox{for $x\in\RR$, $y=0$}\\
u(-\infty,y)=&1,\quad\hbox{uniformly in $y\in\RR_-$}\\
 u(+\infty,y)=&0\quad\hbox{pointwise in $y\in\RR_-$}
\end{array}
\right.
\end{equation}
Note that the convergence of $u$ to 1 to the left, and 0 to the right, are not requested to hold in the same sense. This is not entirely innocent, we will explain why in more detail below. 
We will also be interested in a (seemingly) more complex version of \eqref{e1.2}. We look for a real $c>0$, a function $u(x)$, defined for $x\in\RR$, a function $v(x,y)$ defined in $\RR^2_-:=\RR_+\times\RR_-$, and a curve $\Gamma\subset\RR^2_-$ such that
\begin{equation}
\label{e1.1}
\left\{
\begin{array}{rll}
-d\Delta v+c\partial_x v=&0\quad(x,y)\in\{v>0\}\\
\vert\nabla v\vert=&1\quad((x,y)\in\Gamma:=\partial\{v>0\}\\
\ \\
-Du_{xx}+c\partial_xu+1/\mu u-v=&0\quad\hbox{for $x\in\RR$, $y=0$}\\
v_y=&\mu u-v\quad\hbox{for $x\in\RR$, $y=0$ and $v(x,0)>0$}\\
u(-\infty)=1/\mu,\ u(+\infty)=&0\\
 v(-\infty,y)=&1\quad\hbox{uniformly with respect to $y\leq0$}\\
  v(+\infty,y)=&0\quad\hbox{pointwise with respect to $y\leq0$}.
\end{array}
\right.
\end{equation}
We ask for the global shape of the free boundary $\Gamma$. Before that, we ask about the  existence of a solution $(c,\Gamma,u)$ to \eqref{e1.2}, and of a solution $(c,\Gamma,u,v)$ to \eqref{e1.1}, this indeed deserves some thought, as the condition at $-\infty$ is rather strong.

\noindent Systems \eqref{e1.2} and \eqref{e1.1}  arise from a class of models proposed by H. Berestycki, L. Rossi and the second author to model the speed-up of biological invasions by lines of fast diffusion, see for instance \cite{BCRR} or \cite{BRR2}. The two-dimensional lower half-plane (that was called "`the field"' in the afore-mentionned references), in which reaction-diffusion phenomena occur, interacts with the $x$ axis ("`the road") which has a much faster diffusion $D$ of its own. In Model \eqref{e1.1},  $u(x)$ the density of individuals on the road, and $v(t,x,y)$ the density of individuals in the field. The road yields the fraction $\mu u$ to the field, and retrieves the fraction $\nu v$ in exchange; the converse process occurs for the field.  Model \eqref{e1.2} is obtained from \eqref{e1.1} by forcing $\mu u=v$ on the road, so that the sole unknown is $v(x,y)$, and the exchange term is replaced by $\di\frac{v_y}\mu$ (and $v$ has been renamed $u$). In the sequel, Model \eqref{e1.1} will be called the "model with two species" (that is, the density on the road and in the field may be different), while Model \eqref{e1.2} will be called "Model with one species". Also note that, in both models, the unknown functions are travelling waves of an evolution problem where the term $c\partial_x u$ (resp. $c\partial_xv)$ is replaced by $\partial_tu$ (resp. $\partial_tv$). This is explained in more detail in \cite{CafR},  where the study of \eqref{e1.2} and \eqref{e1.1} was initiated. 
\subsection{Results}
\begin{thm} (Existence for the model with one species) Assume $D\geq d$. 
\label{t1.1}
System \eqref{e1.2} has a solution $(c>0,\Gamma,u)$ with $u$ globally Lipschitz, and we have $\partial_xu<0$. Moreover
\begin{itemize}
\item[--] $\Gamma$ is an analytic curve, as well as a locally Lipschitz graph in the $y$ variable:
\begin{equation}
\label{e1.3}
\Gamma=\{(\varphi(y),y),\ y<0\}.
\end{equation}
\item[--]  $\Gamma\cap\{y=0\}$ is nonempty, assume $(0,0)\in\Gamma$. There is $\e_0>0$ such that
\begin{equation}
\label{e1.4}
\Gamma\cap B_{\e_0}(0)=\{y=\phi(x),\ -\e_0\leq x\leq 0\},\quad \phi(x)=-\frac{x^2}{2D}+o_{x\to0}(x^2).
\end{equation}
\end{itemize}
\end{thm}
\begin{thm} (Existence for the model with two species)
\label{t1.3}
System \eqref{e1.1} has a solution $(c>0,\Gamma,u,v)$ such that the function $v$ is globally Lipschitz, and we have $\partial_xv<0$. Moreover
\begin{itemize}
\item[--] $\Gamma$ is an analytic curve, as well as a locally Lipschitz graph in the $y$ variable:
\begin{equation}
\label{e1.30}
\Gamma=\{(\varphi(y),y),\ y<0\}.
\end{equation}
\item[--] $\Gamma\cap\{y=0\}$ is nonempty, assume $(0,0)\in\Gamma$. There is $\lambda>0$ and $\e_0>0$ such that
\begin{equation}
\label{e1.40}
\Gamma\cap B_{\e_0}(0)=\{y=\phi(x),\ -\e_0\leq x\leq 0\},\quad \phi(x)=\lambda x +o_{x\to0}(\vert x\vert).
\end{equation}
\end{itemize}
\end{thm}

\noindent The main question of this  work, namely, how $\Gamma$ looks like, is addressed in the following theorem. Both models with one species and model with two species are concerned.
Let $c_0$ be the speed of the basic travelling wave $\phi_0(x)$:
\begin{equation}
\label{e2.00}
\begin{array}{rll}
-\phi_0''+c_0\phi_0'=&0\ (x\in\RR_-)\\
\lim_{x\to-\infty}\phi_0(x)=&1,\quad\phi_0(x)=0\ (x\in\RR+)\\
\big[\phi_0\big]_{x=0}=&0,\ \big[\phi_0'\big]_{x=0}=1.
\end{array}
\end{equation}
In what follows, we stress the dependence on $D$ of the velocity, free boundary and solution of the PDE by denoting them $c_C$, $\Gamma_D$,   $u_D$, $v_D$.
\begin{thm}
\label{t1.2} In both models \eqref{e1.2} and \eqref{e1.1}, 
there is $D_0\geq0$ such that, for every $D\geq D_0$ we have $c_D>c_0$. Let $\alpha_D$ be given by $\sin\alpha_D=\di\frac{c_0}{c_D}$. If $\varphi_D$ is given by \eqref{e1.4} we have
\begin{equation}
\label{e1.5}
\lim_{y\to-\infty}{\varphi_D'}(y)=-{\mathrm{cotan}}\alpha_D.
\end{equation}
\end{thm}
In other words, $\Gamma_D$ has an asymptotic direction, which is a line making the angle $\alpha_D$ with the horizontal. We have a more precise version of Theorem \ref{t1.2}:
\begin{thm}
\label{t1.4}
Assume that $\alpha<\di\frac\pi2$.  For every $\omega\in(0, \di\frac12c_D\sin\alpha_D)$, we have
\begin{equation}
\label{e1.50}
{\varphi_D'}(y)=-{\mathrm{cotan}}\alpha_D+O(e^{\omega y}).
\end{equation}
\end{thm}
Thus there is a straight line making the angle $\alpha_D$ with the horizontal that is asymptotic to $\Gamma_D$ at infinity, and the distance between the two shrinks exponentially fast.
\subsection{Underlying question, comments,  organisation of the paper}
\noindent Let us be more specific about the question that we wish to explore here. We want to account for a loss of monotonicity phenomenon for $u$ (Model \eqref{e1.2}) or $v$ (Model \eqref{e1.1}) in the $y$ variable, a study that was initiated in \cite{CafR}. This phenomenon was discovered numerically by A.-C. Coulon in her PhD thesis \cite{ACC}. There, she provided simuations for the evolution system (notice that the propagation takes place upwards):
\begin{equation}
\label{e1.3000}
\begin{array}{rll}
\partial_t u-D \partial_{xx} u= & \nu v(t,x,0)-\mu u \   \    &x\in\R\\
\partial_t v-d\Delta v=&f(v)\    \    &(x,y)\in\R\times\R_+\\
\partial_y v(t,x,0)=&\mu u(t,x,t)-  \nu v(t,x,0)\   \    &x\in\R.
\end{array}
\end{equation}
The   figures below account for some of her results; the parameters are 
$$f(v)=v-v^2,\ \ D=10,\ \ u(0,x)={\bf 1}_{[-1,1]}(x),\ \ v(0,x,y)\equiv0.
$$
 The top figure represents the levels set 0.5 of $v$ at times 10, 20, 30, 40, while the bottom figure represents the shape of $v(40,x,y)$. 

\hfill\includegraphics[width=14cm]{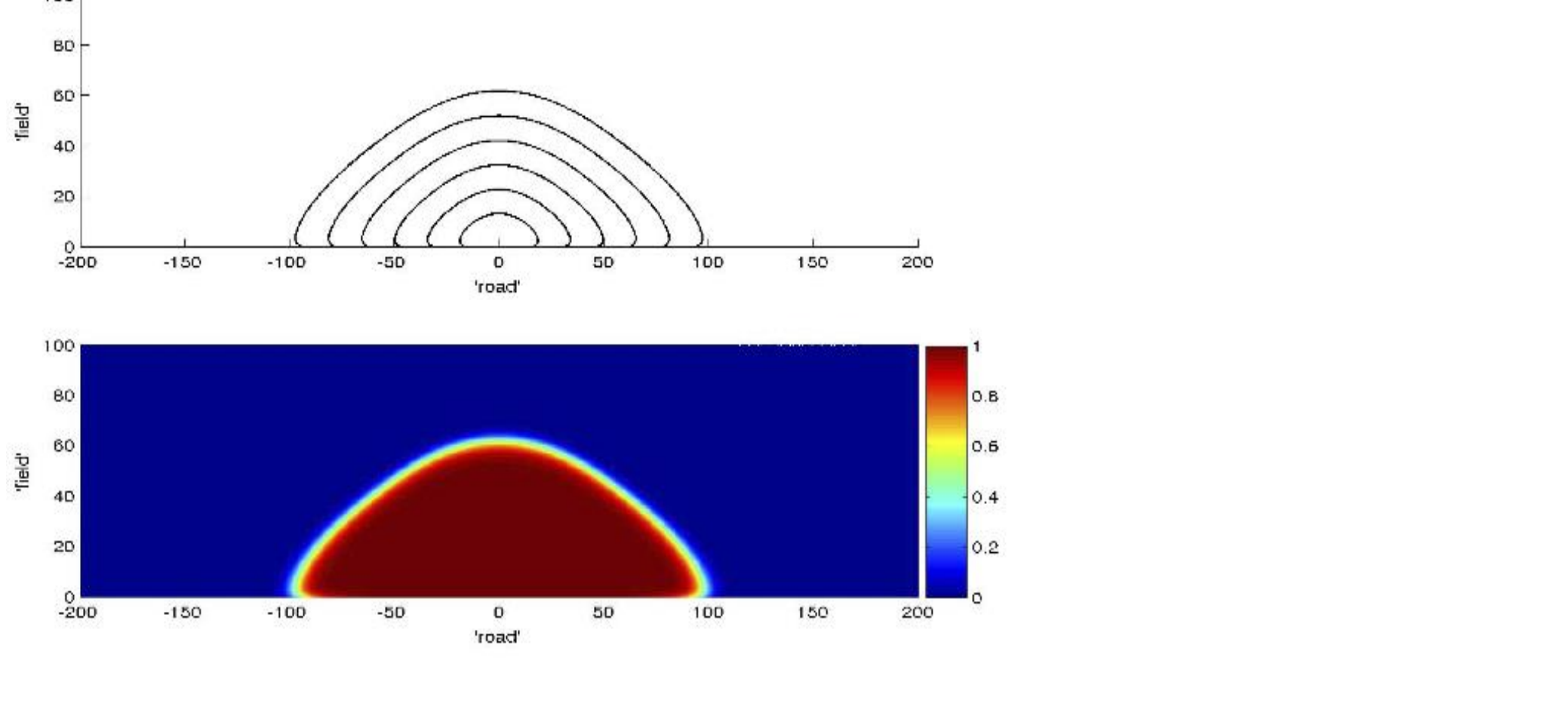} 

\noindent The surprise  is the location of the leading edge of the invasion front: rather than being located on the road, as one would have expected (especially for large $D$), it sits a little further in the field. This entails a counter-intuitive loss of monotonicity. As working directly with the reaction-diffusion \eqref{e1.3000} has not resulted in a significant outcome so far, the idea was to replace the reaction-diffusion model by a free boundary problem, that may be seen as an extreme instance of reaction-diffusion. Using this idea, a first explanation of the location of the leading edge is provided in \cite{CafR}. The conclusions are included in Theorems \ref{t1.1} and \ref{t1.3}.

\medskip
\noindent The goal of the paper is, as already announced, to account for the global shape of the free boundary. We claim that it will provide a good theoretical explanation of A.-C. Coulon's numerics. Indeed it may be expected (although it is not a totally trivial statement) that the solution of \eqref{e1.3000} will converge to a travelling wave. As the underlying situation is that of an invasion, it is reasonable to assume that no individuals (if we refer, as was the initial motivation, to a biological invasion) are present ahead of the front. This is why we impose the seemingly stringent, but reasonable from the modelling point of view, condition at $-\infty$. Theorem \ref{t1.2} shows that it entails a very specific behaviour. 

\noindent Let us explain the consequences of our results on the understanding of the model. Theorems \ref{t1.1} to \ref{t1.2} put together depict a free boundary whose leading edge is in the lower half plane and which, after a possibly nonempty but finite set of turns, becomes asymptotic to a line that goes to the right of the lower half plane. This is in good qualitative agreement with the upper two-thirds of the picture presented here, the lower third accounting for the fact (still to be described in mathematically rigorous terms) that the free boundary bends in order to connect to a front propagating downwards, which is logical as we start from a solution that is compactly supported. Let us, however, point out that the analogy should not be pushed further than what is reasonable, as the logistic nonlinearity displays - and this is especially true for the model \eqref{e1.3000} - more counter-intuitive oddities of its own, see \cite{BRR3}.

\noindent Some additional mathematical comments are in order. The first one is that we have put some effort in proving existence theorems. The reason is that we could not entirely rely on our previous study \cite{CafR} for that: although the leading edge was analysed, we had chosen to wipe out the additional difficulties coming from the study in the whole lower half-plane, by studying a model in a strip of finite length with Neumann conditions at the bottom. While the study of the leading edge is purely local, and will not need any more development, the condtion $\di\lim_{x\to-\infty}u(x,y)=1$ (resp. $\di\lim_{x\to-\infty}u(x,y)=1$) {\it uniformly} in $y\in\RR_-$ requires some additional care, that is presented in Section 3. 

\noindent Theorem \ref{t1.2} is hardly incidental. It is in fact a general feature in reaction-diffusion equations in the plane. The heuristic explanation is the following:   looking down very far in the lower half plane, we may think that the free boundary $\Gamma_D$ propagates like the 1D wave in its normal direction, that is, $V_n=c_0$. On the other hand, it propagates with speed $c_D$ horizontally: this imposes the angle $\alpha_D$.
For a rigorous proof of that, we take the inspiration from previous works on conical-shaped waves for reaction-diffusion equations in the plane. A first systematic study may be found in \cite{BH}, while the stability of these objects is studied in \cite{HMR1}. Further qualitative properties are derived in \cite{HMR3}. Solutions of the one phase free boundary problem for $c\partial_x u-\Delta u=0$ are classified by Hamel and Monneau in \cite{HM}. One of their results will play an important role in the proof of Theorem \ref{t1.2}, this will be explained in detail in Section 4.

\noindent Although this work is clearly aimed at understanding the situation for large $D$ (the case $D\leq d$ poses interesting technical questions in for the one species model) we have not really provided a systematic study of the limit $D\to+\infty$. This will be the object of a forthcoming paper. Another interesting question is whether the free boundary has turning points. While the simulations cleary point at convexity properties of the sub-level sets, we do not, at this stage,   have real hints of what may be true.

\noindent The paper is organised as follows.  In Section 2, we provide some universal bounds on the velocities, in therms of the diffusion on the road $D$. In Section 3, we construct the wave for the one species model and prove Theorem \ref{t1.3}. In Section 4, we indicate the necessary modifications for the two-species model. In Section 5, we prove the exponential convergence of the level sets.
\section{The one species model: Bounds on the velocity in a truncated problem}
 Solutions to \eqref{e1.2} will be constructed through a suitable approximation in a finitely wide cylinder; we set
$$
\Sigma^L=\RR\times(-L,0),
$$
from (a trivial modification of) \cite{CafR}, for $L>1$, there is a solution $(c_D^L,\Gamma_D^L,u_D^L)$ to the auxiliary problem
\begin{equation}
\label{e2.100}
\left\{
\begin{array}{rll}
-d\Delta u+c\partial_x u=&0\quad(x,y)\in\Sigma^L\cap\{u>0\}\\
\vert\nabla u\vert=&1\quad((x,y)\in\Gamma:=\partial\{u>0\}\cap\Sigma^L\\
\ \\
-Du_{xx}+c\partial_xu+1/\mu u_y=&0\quad\hbox{for $x\in\RR$, $y=0$}\\
u(-\infty,y)=&1,\quad u(+\infty,y)=(1-y-L)^+\\
u(x,-L)=&1.
\end{array}
\right.
\end{equation}
In biological terms, this means that the boundary $\{y=-L\}$ is lethal for the individuals. The limit at $x=+\infty$, namely the function $(1-y-L)^+$, is of course not chosen at random, it solves the one-dimensional free boundary problem
$$
\left\{
\begin{array}{rll}
-u''=&0\quad\hbox{in $(-L,0)\cap\{u>0\}$}\\
\vert u'\vert=&1\quad\hbox{on $\partial\biggl((-L,0)\cap\{u>0\}$}\biggl)\ \hbox{(that is, at $y=1-L$)}\\
u(-L)=&0.
\end{array}
\right.
$$
To ensure the maximum chance to retrieve, in the end, a solution that converges to 1 uniformly in $y$ as $x\to-\infty$, we have imposed the Dirichlet condition $u=1$ at the bottom of the cylinder.

\subsection{Exponential solutions}
\noindent At this point, we need to make a short recollection of what the exponential solutions of the linear problem are. System \eqref{e2.100}, linearised around 0, reads
\begin{equation}
\label{e2.60}
\left\{
\begin{array}{rll}
-d\Delta u+c\partial_x u=&0\quad(x,y)\in\Sigma^L\\
\ \\
-Du_{xx}+c\partial_xu+1/\mu u_y=&0\quad\hbox{for $x\in\RR$, $y=0$}\\
u(x,-L)=&0\\
u(-\infty,y)=&0.
\end{array}
\right.
\end{equation}
It has exponential solutions that decay to 0 as $x\to-\infty$, i.e. solutions of the form $\psi_D^L(x,y)=e^{\alpha x}h(y)$ (we have to put here the dependence with $L$ and $D$. If $L=+\infty$ we may choose the function $h$ as an exponential in $y$.  We have
$$
\psi^L_D(x,y)=e^{\alpha^L_D x}\mathrm{sh}(\beta^L_D(y+L)),
$$
so that the exponents $\alpha^L_D$ and $\beta^L_D$ satisfy
\begin{equation}
\label{e2.3001}
\left\{
\begin{array}{rll}
-d(\alpha^2+\beta^2)+c\alpha=&0\\
-D\alpha^2+c\alpha+\di\frac\beta\mu{\mathrm{cotanh}}(\beta L)=&0.
\end{array}
\right.
\end{equation} 
Three types of limits will be considered.

\noindent{\bf Case 1.} The limit $D\gg L\gg1$, $c$ bounded. We expect $\beta_D^L$ to go to 0 as $D\to+\infty$, so that 
$$
\beta_D^L{\mathrm{cotanh}}(\beta_D^LL)\sim\frac1L.
$$
Then \eqref{e2.3001} yields 
estimates of the form 
\begin{equation}
\label{e2.7}
\alpha^L_D\sim\sqrt{\frac1{\mu LD}},\quad \beta^L_D\sim\sqrt{\frac{c}{d\sqrt{\mu DL}}}.
\end{equation}
It is to be noted that these equivalents may be pushed up to $c=o(\sqrt D)$.  

\noindent{\bf Case 2.}  The limit $D\gg L\gg1$, $c\gg\sqrt D$. This time we  expect $\beta_D^L$ to go to infinity as $D\to+\infty$, so that 
$$
{\mathrm{cotanh}}(\beta_D^L)\to 1.
$$
We have
\begin{equation}
\label{e2.8}
\alpha_D^L\sim\biggl(\frac{c}{dD^2}\biggl)^{1/3},\quad \beta_D^L\sim\mu\biggl(\frac{c}{\sqrt D}\biggl)^{2/3}.
\end{equation}
We have here a first occurrence of the critical order of magnitude $c\propto\sqrt D$.

\noindent{\bf Case 3.}  The limit $L\gg D\gg1$, $c\gg1$. We expect that $L\beta_D^L$ will go to infinity, so that ${\mathrm{cotanh}}(\beta_D^L)\to 1.$ In this setting, we have the estimates \eqref{e2.8}.

\subsection{Universal upper bound}
In the construction of a travelling wave for \eqref{e1.2}, the first task is to bound the velocity from above. We will prove straight away the upper bound that will serve us in the later section, namely that $c_D^L$ cannot exceed a (possibly large) multiple of $\sqrt D$. 
\begin{thm}
\label{t2.10}
There is $K>0$, independent of $D$ and $L$, such that 
\begin{equation}
\label{e2.9}
c_D^L\leq K\sqrt D.
\end{equation}
\end{thm}
\noindent{\bf Proof.} Assume $c_D^L\gg\sqrt D\gg L\gg1$, so that \eqref{e2.8} hold for $\alpha_D^L$ and $\beta_D^L$. Set
$$
\underline u_D^L(x,y)=1-\psi_D^{2L}(x,y).
$$
The function $\underline u_D^L$ vanishes on the curve $\underline\Gamma_D^{2L}$ whose equation is
$$
x=-\frac1{\alpha_D^{2L}}{\mathrm{ln}}\ {\mathrm{sh}}\beta_D^{2L}(y+2L).
$$
In particular, for $y=-L$ we have
\begin{equation}
\label{e2.10}
x=-\frac1{\alpha_D^{2L}}{\mathrm{ln}}\ {\mathrm{sh}}(\beta_D^{2L}L).
\end{equation}
Translate the solution $u_D^L(x,y)$ of \eqref{e1.2} so that the leftmost point of $\Gamma_D^L$ is located at the vertical of the origin. Then, from \eqref{e2.10}, $\underline\Gamma_D^{2L}$ lies entirely to the left of $\Gamma_D^L$. By the maximum principle we have
$$
u_D^L(x,y)\geq\underline u_D^L(x,y).
$$
Now, translate $\underline u_D^L$ as much as this allows it to remain under $u_D^L$. If $\lambda_D^L$ is the maximum amount by which one may translate $\underline u_D^L$,  there is a contact point $(x_D^L,y_D^L)$  between $u_D^L(x,y)$ and $\underline u_D^L(x-\lambda_D^L,y)$. 

\noindent Still from \eqref{e2.10}, we have $y_D^L>-L$. So, it now depends on whether we have $y_D^L=0$ or $y_D^L\in(-L,0)$. In the first case,  let $\Gamma_D^L$ intersect the line $\{y=0\}$ at the point $(\tilde x_D^L,0)$. If 
$x_D^L=\tilde x_D^L$, we have $\partial_x u_D^L(x_D^L,0)=0$ but we also have
$$
\partial_x u_D^L(x_D^L,0)<\partial_x\underline u_D^L(x_D^L,0)<0,
$$
an impossibility. If not, we have $x_D^L<\tilde x_D^L$ and this time we have
$$
\partial_x u_D^L(x_D^L,0)=\partial_x\underline u_D^L(x_D^L,0),\quad \partial_{xx} u_D^L(x_D^L,0)\geq\partial_x\underline u_D^L(x_D^L,0),\quad \partial_y u_D^L(x_D^L,0)<\partial_y\underline u_D^L(x_D^L,0),
$$
the last inequality because of the Hopf lemma. This prevents the Wentzell relation from holding at $(x_D^L,0)$. So, the only possibility left is $y_D^L<0$. By the strong maximum principle, the point $(x_D^L,y_D^L)$ is on the free boundary $\Gamma_D^L$. We have
$$
\partial_y u_D^L(x_D^L,y_D^L)=-\beta_D^{2L}e^{\alpha_D^{2L}(x_D^L-\lambda_D^L)}{\mathrm{ch}}\beta_D^{2L}(y_D^L+2L),
$$
that is,
$$
\vert\partial_y u_D^L(x_D^L,y_D^L)\vert\geq\beta_D^{2L}e^{\alpha_D^{2L}(x_D^L-\lambda_D^L)}{\mathrm{sh}}\beta_D^{2L}(y_D^L+2L)=\beta_D^{2L},
$$
because $\underline u_D^L(x_D^L,y_D^L)=0$. Thus, if $D$ is sufficiently large, 
\eqref{e2.8} implies 
$$\vert \nabla \underline u_D^L(x_D^L,y_D^L)\vert>1,
$$
contradicting the free boundary relation for $u_D^L$. This proves the theorem. \hfill$\Box$
\subsection{Universal lower bound}
\begin{thm}
\label{t2.1}
There is $L_0>0$ such that the velocity $c_D^L$ in Model \eqref{e2.100} satisfies
$$
\lim_{D\to+\infty}c^L_D=+\infty,
$$
uniformly with respect to $L\geq L_0$.
\end{thm}

\noindent The proof of Theorem \ref{t2.1} will be by contradiction. From now on, and until this has been proved wrong, we assume that, $c^L_D$ is bounded, both with respect to $L$ and $D$.   Recall that $\Gamma^L_D$, the free boundary, is an analytic  graph $\{x=\varphi^L_D(y)\}$. Moreover, from \cite{CafR}, we may always assume that it intersects the line $\{x=0\}$, so that we may always assume $\varphi^L_D(0)=0$. Define $\bar x^L_D$ as the last $x$ such that, for all $y\in(-L,0)$, then $(x,y)\notin\Gamma^L_D$.   Our main step is to prove that the front goes far to the left of the domain, this is expressed by the following lemma.
\begin{lm}
\label{l2.1}
There is $q>0$ universal such that
\begin{equation}
\label{e2.15}
\bar x^L_D\leq -q\sqrt D.
\end{equation}
\end{lm}
\noindent{\sc Proof.} Recall that we have 
$$
\di\frac{d\varphi^L_D}{dy}(0)>0.
$$
For every $x<0$, let $y^L_D(x)$ be the first $y$ such that $x=\varphi^L_D(y)$. We will prove that
\begin{equation}
\label{e2.2}
\lim_{\e\to0}y^L_D(-\e\sqrt D)=0,\ \hbox{uniformly in $D\geq2$ and $L\geq2$},
\end{equation}
which implies Lemma \ref{l2.1}.
Assume \eqref{e2.2} does not hold, and consider $\delta_0>0$ such that, for a sequence $(D_n,L_n)_n$, going to $+\infty$, and a sequence $(\e_{n})_{n}$
going to $0$ as $n\to+\infty$, we have
$$
y^{L_n}_{D_n}(\e_{n}\sqrt{D_n})\leq-\delta_0.
$$
Obviously, we must assume the boundedness from below of the sequence $(\e_n\sqrt {D_n})_n$. For every $n$ set 
$$
X_n=(-\e_{n}\sqrt{D_{n}},-\delta_0/2),\ \bar X_n=(-\e_{n}\sqrt{D_{n}},0).
$$
 For $\e_0$, to be chosen small in due time, the segment $(\bar X_n,X_{n})$ is at distance at least $\delta_0/2$ from the free boundary. By nondegeneracy (recall the boundedness of $(c_{D_n}^{L_n})_n$), there is $q_0>0$ universal such that, for all $n$: 
\begin{equation}
\label{e2.3}
u_{D_n}^{L_n}(X)\geq q_0\ \hbox{on $[\bar X_n,X_n]$}.
\end{equation}

\noindent On the other hand, recall that $\partial_yu_{D_n}^{L_n}(x,0)\leq1$, this follows from \cite{CafR}. So, the equation for $u_{D_n}^{L_n}(x,0)$ reads, simply
$$
\partial_{xx}u_{D_n}^{L_n}(x,0)\leq O(D^{-1}),\ u(0,0)=u_x(0,0)=0.
$$
Thus we have 
$$
u_{D_n}^{L_n}(-\e_n\sqrt D_{n},0)\leq O(\e_n).
$$
This yields, for $n$ large enough, the existence of $\delta_0''>0$, universal, such that 
\begin{equation}
\label{e2.4}
\partial_yu_{D_n}^{L_n}(x,0)\leq-\delta_0',\ x\in(-\e_0\sqrt{D_{n}},\e_{n}\sqrt{D_{n}}).
\end{equation}
This comes from the Hopf boundary lemma. So now, we now write the equation for $u_{D_n}^{L_n}(x,0)$ as 
$$
\partial_{xx}u_{D_n}^{L_n}(x,0)-\frac{o(1)}{\sqrt D}\partial_xu_{D_n}^{L_n}(x,0)=\frac{\partial_yu_{D_n}^{L_n}(x,0)}D\leq-\frac{\delta_0}D,\quad u_{D_n}^{L_n}(0,0)=\partial_xu_{D_n}^{L_n}(0,0)=0.
$$
Integrating this equation on $(-\e_0\sqrt{D_n},0)$ and invoking \eqref{e2.4} allows us to find a small constant $\delta_0''>0$ such that
 $$u_{D_n}^{L_n}(-\e_0\sqrt {D_{n}},0)\leq-\delta_0'',
$$
 a contradiction. \hfill$\Box$

\noindent{\sc Proof of Theorem \ref{t2.1}.} Recall that we have still assumed the boundedness of $(c^{L_n}_{D_n})_{n}$. In order to allevite the notations a little, we omit the index ${}_n$. Integration of \eqref{e2.100}
over $\Sigma^L$ yields
 \begin{equation}
\label{e2.5}
(1/\mu+L)c^L_D=\int_{\Gamma^L_D}\partial_\nu u_D^L-\int_{-\infty}^{+\infty}\vert \partial_yu_D^L(x,-L)\vert dx.
\end{equation}
This expression has to be handled with care, because each integral, taken separately, diverges. We will see, however, that there is much less nonsense in \eqref{e2.5} than it carries at first sight. The curve $\Gamma^L_D$ has an upper branch, that we call $\Gamma^L_{D,+}$, and that connects $(0,0)$ to $X^L_D$; the latter point being a turning point. The lower branch, called $\Gamma^L_{D,-}$, connects $X^L_D$ to $(+\infty,-L+1)$, in other words it is asymptotic to the line $\{y=-L+1\}$ as $x$ goes to $+\infty$. We decompose $\eqref{e2.5}$ into 
$$
\begin{array}{rll}
&(1/\mu+L)c^L_D\\
=&\di{\int_{\Gamma^L_{D,+}}\partial_\nu u_D^L-\int_{-\infty}^{x^L_D}\vert \partial_yu_D^L(x,-L)\vert dx}+\biggl(\int_{\Gamma^L_{D,+}}\partial_\nu u_D^L-\int_{x^L_D}^{+\infty}\vert \partial_yu_D^L(x,-L)\vert dx\biggl)\\
:=&I-I\!I+I\!I\!I.
\end{array}
$$
Let us remark that $\partial_yu_D^L(x,-L)\geq-1$ for all
$x\in\RR$. Indeed, the function $u_D^L(x,y)$ being decreasing in $x$, it is larger than $(1-L-y)^+$. So, we have $u_D^L(x,-L)\geq-1$. Now, the graph $\Gamma^L_D$ has a discrete set of turning points, due to its analyticity. Away from these turning points, for $x\geq x^L_D$, there is a finite number of $y's$: $y_1>y_2>...y_{n_D(x)}$ such that $(x,y_i)\in\Gamma^L_D$ for $1\leq i\leq n_D$. The point $(x,y_1)$ has already been counted in the integration over $\Gamma^L_{D,+}$, and so does not need to be counted again. Notice then that, if $x$ is not the abscissa of a turning point, then $n_D$ (or $n_D-2$) is even, because of the configuration of $\Gamma^L_D$. So, $u_\nu(x,y_i)=u_\nu(x,y_{i+1})=1$ for $1\leq i\leq n_D(x)-2$: it  suffices to consider the sole point $(x,y_{n_D(x)})$ in the computation of $I\!I\!I$, however at this point we also have $\partial_\nu u_D^L=1$. Therefore we end up with $I\!I\!I\geq0$. Of course the integral giving $I\!I\!I$ converges, but we do not even have to bother to prove it.

\noindent From Lemma \ref{l2.1}, we have
$$
I\geq C\sqrt D,
$$
for some universal $C>0$. We already saw that $I\!I\!I$ was nonnegative, so let us deal with $I\!I$. By nondegeneracy, there is $A>0$ universal such that 
$$
u_D^L(x,y)\geq 1-A\psi_D^L(x,y),\quad x\leq x^L_D-1.
$$
Notice that we do not change anything if, in $I\!I$, we integrate up to $x^L_D-1$ instead of $x^L_D$. We deduce, because $\psi_D^L(x,-L)=0$ and 
$u_D^L(x,-L)=1$: 
$$
\partial_yu_D^L(x,-L)\geq -A\partial_y\psi_D^L(x,-L),\quad x\leq x^L_D-1.
$$
This implies
$$
I\!I\leq O(1)+O({c}{\sqrt{DL}})^{1/2}.
$$
This yields
$$
\frac{c^L_D}{\sqrt{DL}}\geq \frac{q}{\sqrt L}-O(\frac{c}{\sqrt {DL}})^{1/2},
$$
for a universal $q$, as soon as we choose $D\gg L\gg1$. This is an obvious contradiction.  Now, note that, because of the Dirichlet condition at $y=-L$, the sequence $(c^L_D)^L$, for fixed $D$,  is increasing.  This ends the proof of the theorem.\hfill$\Box$

\section{The one species model: construction and properties of the free boundary}
\subsection{Global solutions of the free boundary problem in the plane}
In this short section we recall a result of Hamel and Monneau that we will use to analyse the behaviour of the free boundary at infinity. Consider a solution $(c>0,\Gamma,u)$ of the free boundary problem in the whole plane
\begin{equation}
\label{e4.1}
\begin{array}{rll}
-\Delta u+c\partial_xu=&0\quad\hbox{in $\Omega:=\{u>0\}$}\\
\vert\nabla u\vert=&1\quad\hbox{on $\Gamma:=\partial\Omega$}.
\end{array}
\end{equation}
The result is a classification of the solutions of \eqref{e4.1} having certain additional properties. We rephrase it here to avoid any confusion, since the function $u$ in \cite{HM} corresponds to $1-u$ in our notations.
\begin{thm} (Hamel-Monneau \cite{HM}, Theorem 1.6)
\label{t4.1}
Assume that
\begin{enumerate}
\item $\Gamma$ is a $C^{1,1}$ curve with globally bounded curvature,
\item $\RR^2\backslash\Omega$ has no bounded connected component,
\item we have 
\begin{equation}
\label{e4.2}
\liminf_{d(X,\Gamma)\to+\infty, X\in\Omega}u(X)=1.
\end{equation}
\end{enumerate}
Then $c\geq c_0$ and, if we set
\begin{equation}
\label{e4.4}
 \sin\alpha=\frac{c_0}c,
 \end{equation} 
 either $u$ is the tilted one-dimensional solution $\phi_0(y\cos\alpha\pm x\sin\alpha)$, or $u$ is a conical front with angle $\alpha$ to the horizontal, that is, the unique solution $u_\alpha(x,y)$ of \eqref{e4.1} such that $\Gamma$ is asymptotic to the cone
 $\partial\mathcal{C_\alpha}$, with
\begin{equation}
\label{e4.9}
{\mathcal{C}}_\alpha=\{(x,y)\in\RR_-\times\RR:\  \frac{x}{\vert y\vert}\leq-{\mathrm{cotan}}\alpha\}.
\end{equation}
\end{thm}
Note that the fact that $u_\alpha$ is unique is not exactly trivial, it is given by Theorem 1.3  of \cite{HM}. Let us already notice that Properties 1 and 2 of this theorem are satisfied by the solution $(c_D^L,\Gamma_D^L,u_D^L)$ of \eqref{e2.100}: Property 1 is clear, and Property 2 is readily granted by the monotonicity in $x$. Property 3 will be a little more involved to check.
\subsection{Construction of a solution in the whole half-plane}
From then on, we fix $D>0$ large enough so that Theorem \ref{t2.1} holds. Pick $L_0$ such that $c_D^{L_0}>c_0$, this implies $c_D^L\geq c_D^{L_0}>c_0$. 
Notice that we have almost all the elements for the proof of Theorem \ref{t1.1}, we just need, in addition, to control where the free boundary meets the line $\{y=0\}$.  Since we have this freedom, we assume $\Gamma_D^L$ to intersect the line $\{y=0\}$ at the origin. Let $\bar X_D^L=(\bar x_D^L,\bar y_D^L)$ the point of $\Gamma_D^L$ that is furthest to the left, our sole real task will be to prove that $X_D^L$ cannot escape too far as $L\to+\infty$. Indeed, we notice that the property
\begin{equation}
\lim_{x\to-\infty}u_D^L(x,y)=1,\quad\hbox{uniformly with respect to $y<0$ and $L\geq L_0$}
\end{equation}
holds easily. Indeed, fro the maximum principle we have
\begin{equation}
\label{e4.3}
1\geq u_D^L(x,y)\geq(1-e^{\gamma x})^+,
\end{equation}
for all $\gamma<\di\frac{c_D^{L_0}}D$, as it is a subsolution to the equation for $u$ in the plane and on the line $\{y=0\}$, and below $u_D^L$ on the bottom line and also the vertical segment $\{x=0,-L\leq y\leq0\}$. Hence, at that point, we have almost everything for the construction of the wave in the whole half-plane, except the attachment property.
\begin{lm}
\label{l4.1}
There is a constant $K_D>0$ independent of $L$ such that $\vert \bar X_D^L\vert\leq K_D$.
\end{lm}
\noindent{\bf Proof.} Let us first assume that
$$
\lim_{L\to+\infty}\bar y_D^L=-\infty.
$$
Translate $\Gamma_D^L$ and $u_D^L$ so that $\bar X_D^L$ becomes the new origin, the free boundary $\Gamma_D^L$ meets therefore the horizontal line at the point $(-\bar x_D^L,-\bar y_D^L)$. Up to a subsequence the triple $(c_D^L,\Gamma_D^L,u_D^L)$ converges to a solution $(c_D^\infty,\Gamma_D^\infty,u_D^\infty)$ of the free boundary problem \eqref{e4.1} in the whole plane. Moreover, the origin is its leftmost point. So, we may slide $\phi_0$ from $x=-\infty$ to the first point where it touches $u_D^\infty$, this can only be at $\Gamma_D^\infty$, thus at the origin. But then we have 
$$
\phi_0'(0)=\partial_xu_D^\infty(0,0)=-1,
$$
a contradiction with the Hopf Lemma. So this scenario is impossible, and the family $(\bar y_D^L)_{L\geq L_0}$ is bounded as $L\to\infty$. 

\noindent To prove that the family $(\bar x_D^L)_{L\geq L_0}$ is bounded, we consider
$$
\bar L\geq 2\limsup_{L\to+\infty}(-\bar y_D^L),
$$
and integrate the equation for $u_D^L$ on $\Sigma_{\bar L}=\RR\times(-\bar L,0)$, we obtain
$$
(\frac1\mu+\frac1{\bar L})c_D^L=\int_{\Gamma_D^L\cap\Sigma_{\bar L}}\partial_\nu u_D^L+\int_{\{u(x,-\bar L)>0\}}\partial_yu_D^L(x,-\bar L)dx:=I+I\!I.
$$
From \eqref{e4.3} and elliptic estimates, we have
$$\vert\partial_yu_D^L(x,-\bar L)\vert\leq Ce^{\gamma x},
$$
moreover the choice of $\bar L$ implies that the rightmost point of $\{u(x,-\bar L)>0$ is at the left of the origin. Hence $I\!I$ is uniformly bounded with respect to $L$. On the other hand we have
$$
I\geq -\bar x_D^L,
$$
which, from the  uniform boundedness of $c_D^L$ in Theorem \ref{t2.10}, yields the boundedness of the family $(\bar x_D^L)_{L\geq L_0}$. \hfill$\Box$

\noindent{\bf Proof of Theorem \ref{t1.1}.} We send $L$ to infinity, a sequence $(c_D^{L_n},\Gamma_D^{L_n},u_D^{L_n})_n$ will converge to a solution $(c_D,\Gamma_D,u_D)$ of \eqref{e1.2}. Because of Lemma \ref{l4.1}, the free boundary $\Gamma_D$ meets the line $\{y=0\}$ at a point $(-\bar x_D,0)$ and the expansion \eqref{e1.4} is granted by Theorem 1.4 of \cite{CafR}. Notice also that the uniform limit at $-\infty$ is also granted because of \eqref{e4.3}. So, to finish the proof of the theorem, it remains to prove that, for all $y<0$, the positivity set of $u$ only extends to a finite range. Such were it not the case, the limit
$$
u_D^\infty(x)=\lim_{x\to+\infty}u_D(x,y)
$$
would exist and be nonzero. It would solve the free boundary problem
$$
-(u_D^\infty)''=0\ \hbox{on the positivity set},\quad \vert (u_D^\infty)'\vert=1\ \hbox {at the free boundary points}.
$$
This only allows for   $u^\infty_D(y)=(y+a)^-$, a contradiction to the bounedness of $u_D^\infty$. \hfill$\Box$
\subsection{The tail at infinity}
Let $(c_D,\Gamma_D,u_D)$ the solution constructed in Theorem \ref{t1.1} as one of the limits $L\to+\infty$ of $(c_D^L,\Gamma_D^L,u_D^L)$, with
$$
\Gamma_D=\{(\varphi_D(y),y),\ y\leq0\},
$$
where $\varphi_D$ is a smooth, locally Lipschitz function. 
Theorem \ref{t2.1} readily implies the first part of Theorem \ref{t1.2}, that is
$$
\lim_{D\to+\infty}c_D=+\infty,
$$
simply because $c^L_D\leq c_D$. The only use of this result that we are going to make in this section is that  there exists $D_0>0$ such that $c_D>c_0$ for all $D\geq D_0$, which is the first part of Theorem \ref{t1.2}. And so, we drop the indexes ${}_D$ in the rest of this section, for the simple reason that the dependence with respect to $D$ will not appear anymore. To prove the second part of Theorem \ref{t1.2}, we apply Theorem \ref{t4.1} to any sequence of translates of $u$:
\begin{equation}
\label{e4.5}
u_n(x,y)=u(\varphi(y_n)+x,y_n+y),
\end{equation}
to infer  that any possible limit of $u_n$ is the one-dimensional wave, tilted in the correct direction. Properties 1 and 2 of the theorem being readily true, we concentrate on Property 3. The main step will be to prove that $\phi$ is in fact globally Lipschitz in (any sub-plane of) the half plane, once this is done a suitably designed Hamel-Monneau type \cite{HM} subsolution, placed under $u$, will give the property.
\begin{prop}
\label{p4.1}
The function $\varphi$ is globally Lipschitz in $(-\infty,y_0]$, for any $y_0<0$.
\end{prop}
\noindent{\bf Proof.} Note that the lemma is trivially false if we inisist in making $y$ vary on the whole half-line. Also, as will be clear from the proof, the value of $y_0$ will play no role as soon as it remains a little away from 0. So, we will assume for definiteness $y_0=-1$. The main step of the proposition consists in proving that no point of $\Gamma\cap\{y\leq-1\}$ may have a horizontal tangent. Assume that there is such a point $X_0=(\varphi(y_0),0)$, $y_0\leq -1$. Translate $u$ and $\Gamma$ so that it becomes the origin, still denoting them by $\Gamma$ and $u$. In a neighbourhood of $0$ of size, say, $\rho>0$, $\Gamma$ may be written (recall that it is an analytic curve) as 
$$
\Gamma\cap B_\rho=\{(x,\psi(x),\vert x\vert\leq\rho\}.
$$
We have $\psi(0)=\psi'(0)=0$. Two cases have to be distinguished. The first one is $\psi\equiv0$ in $[-\rho,\rho]$. By analyticity and Cauchy-Kovalevskaya's Theorem, this implies that $\psi$ is defined and equal to 0 on the whole line; actuallly this case may happen only if $X_0$ is a point at infinity, that is, $u$ is a limit of translations of infinite size. But then we have
$$
u(x,y)=y^+\ \hbox{or}\ u(x,y)=y^-,
$$
depending on whether the positivity set of $u$ is above or below $\Gamma$. This is in contradiction with the boundedness of $u$. The other case is $\psi$ nonconstant in $[-\rho,\rho]$, so that $u$ is nonconstant either in $B_\rho$. By analyticity again, $u$ has an expansion of the following type, in a neighbourhood of the origin:
\begin{equation}
\label{e4.6}
u(x,y)=x^n+yP(x)+y^2R(x,y),
\end{equation}
with $n\geq 2$, the functions $P$ and $R$ being smooth in their arguments.  Because $\partial_xu$ it maximal at the origin, we have $u_{xy}(0,0)\neq0$ from the Hopf Lemma. Consider the situation where we have, for instance
$
u_y(0,0)=-1,$ $u_{xy}(0,0)=\beta>0.
$
We have therefore 
$$
P(x)=-1+\beta x+O(x^2),
$$
and
$$
\psi (x)=
\frac{x^n}{-1+\beta x+O(x^2)}.
$$
Inside $B_\rho$ we also have
$$
\nabla u(x,y)=\left(\begin{array}{rll}
&nx^{n-1}+y(2\beta+O(x)+y\partial_xQ(x,y)\\
&-1+2\beta x+O(x^2)+\partial_y(y^2R)(x,y)
\end{array}
\right)
$$
so that
$$
\nabla u(x,\psi(x)))=\left(\begin{array}{rll}
&nx^{n-1}+O(x^n) \\
&-1+2\beta x+O(x^2)
\end{array}
\right)
$$
This implies, because $n\geq 2$:
$$
\vert\nabla u(x,\psi(x))\vert^2=(1-2\beta x+O(x^2))^2+O(x^{2n-2})=1-4\beta x+O(x^2).
$$
Thus the free boundary relation cannot be satisfied in $B_\rho$, except at the origin.

\noindent We note that these situations exhaust what can happen on $\Gamma$. Indeed, if there is a sequence $(X_n)_n$ of $\Gamma$, going to infinity, such that the tangent to $\Gamma$ at $X_n$ makes an angle $\alpha_n$ with the horizontal, with $\di\lim_{n\to+\infty}\alpha_n=0$, the usual translation and compactness argument yields a solution of  the free boundary problem \eqref{e4.1}, where the tangent at the origin is horizontal. Once again we are in one of he above two cases, that are impossible. 
\hfill$\Box$

\noindent The last step is to check a property that will imply Property 3 for any limiting translation of $(\Gamma,u)$ of the form \eqref{e4.5}. This is the goal of the next proposition.
\begin{prop}
\label{p4.2}
.We have
\begin{equation}
\label{e4.7}
\liminf_{u(x,y)>0,y\to-\infty,d((x,y),\Gamma)\to+\infty}u(x,y)=1.
\end{equation}
 \end{prop}
 \noindent{\bf Proof.}  We use the notations of Theorem \ref{t4.1}. Assume for definiteness that the leftmost point of $\Gamma$ is located at $X_l=(0,y_l)$ with $y_l<0$. Pick $\bar X\in\Omega$, call $M>0$ its distance to $\Gamma$, and also set $X_0=(x_0,-N)$,  both $M$ and $N$ will be assumed to be large, independently of one another.  
From Proposition \ref{p4.1} we claim the existence of a cone ${\mathcal{C}}_\beta$ (the notation is given by \eqref{e4.9}), the angle $\beta>0$ depending on the Lipschitz constant of $\varphi$, but independent of $M$ and $N$, and a point $\bar X_{M,N}\in\Omega$,
 such that
 \begin{enumerate}
 \item we have $d(\bar X_{M,N},\Gamma)\geq\sqrt{\inf(M,N)}$,
 \item we have $\bar X\in \bar X_{M,N}+{\mathcal{C}}_\beta$,
 \item the upper branch of $ \bar X_{M,N}+\partial{\mathcal{C}}_\beta$ meets the line of fast diffusion $\{y=0\}$ at a point $(x_{M,N},0)$ with
 \begin{equation}
 \label{e4.10}
 x_{M,n}\leq-\gamma\inf(M,N),
 \end{equation}
 $\gamma>0$ independent of $M$ and $N$.
  \end{enumerate}
  \begin{figure}[htbp] 
    \centering
    \includegraphics[width=15cm]{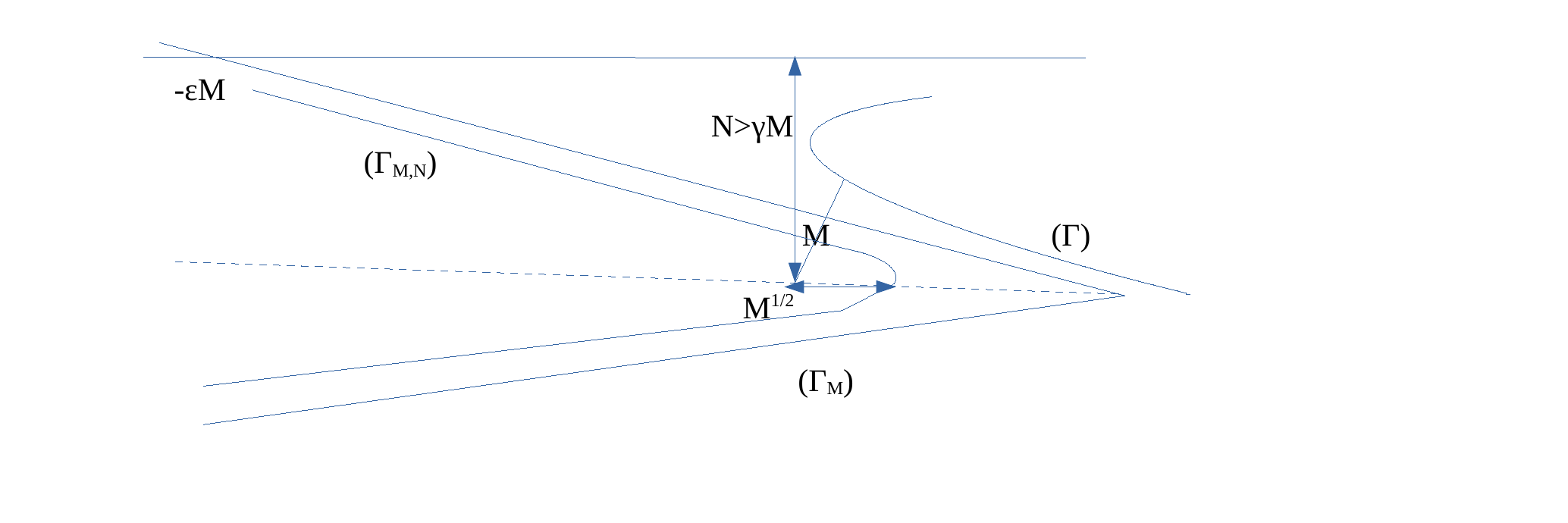} 
    \caption{$\bar X$ below $\Gamma$}
    \label{f4.1}
 \end{figure}
 \begin{figure}[htbp] 
    \centering
    \includegraphics[width=15cm]{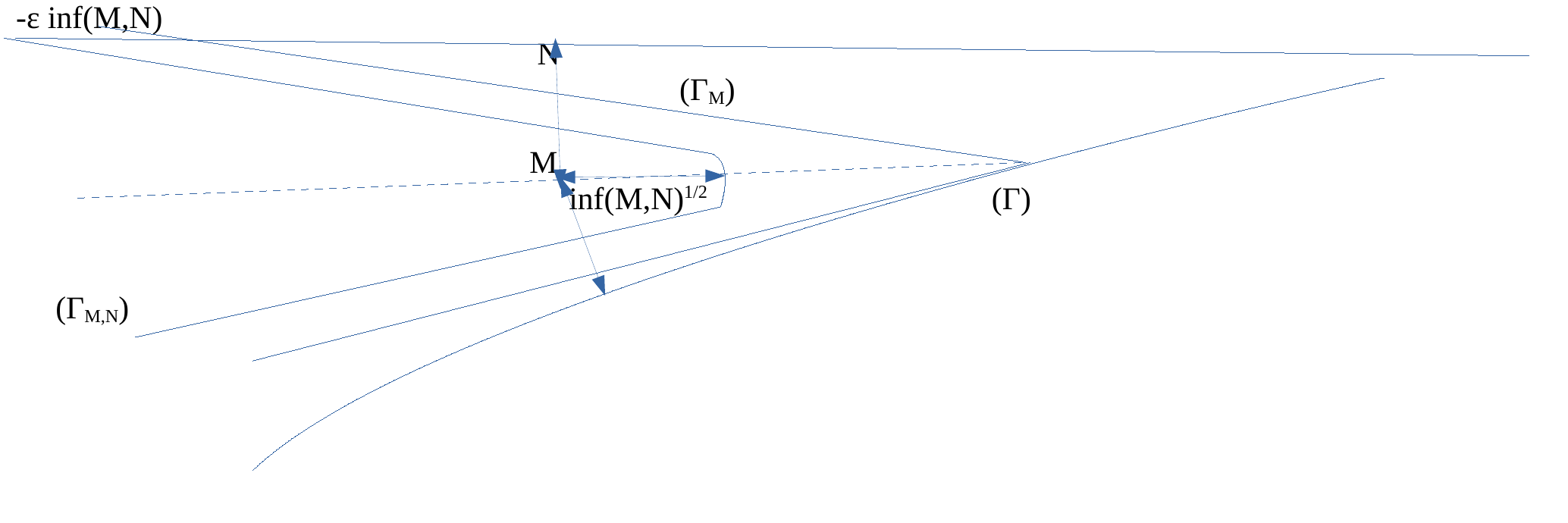} 
    \caption{$\bar X$ above $\Gamma$}
    \label{f4.2}
 \end{figure}
  Two cases have to be distinguished, depending on whether $\bar X$ is below (see Fig. \ref{f4.1}) or above (see Fig. \ref{f4.2}) its projection onto $\Gamma$ (or, if the projection is not unique, the projection that we have selected). Note that, in the first case, there is $C>0$, once again independent of $M$ and $N$, such that 
  $$
  N\geq CM,
  $$
  in such a case we have $\inf(M,N)=M$. We now translate the origin to $\bar X_{M,N}$. Note that the line of fast diffusion becomes the line $\{y=N\}$.
  Consider a smooth, nonpositive, concave, even function $\varphi_{M,N}(y)$ such that 
  \begin{enumerate}
  \item we have $\varphi_{M,N}(0)=0$,
  \item we have $-\varphi_{M,N}''(y)\leq \di\frac1{\inf(M,N)}$,
  \item we have $\vert\varphi_{M,N}'(y)\vert=2\mathrm{cotan}\beta$ if $\vert y\vert\geq\sqrt{\inf(M,N)}$.
  \end{enumerate}
  Let $\Gamma_{M,N}$ be the graph
  $$
  \Gamma_{M,N}=\{(x,y)\in\RR_-\times\RR: x=-\sqrt{\inf(M,N)}+\varphi_{M,N}(y)\}.
  $$
It meets the line of fast diffusion at a point $\tilde x_{M,N}$ that satisfies an estimate of the form
\begin{equation}
\label{e4.11}
\tilde x_{M,N}-x_{M,N}\leq\gamma\inf(M,N),
\end{equation}
with a universal $\gamma>0$, possibly different from that of \eqref{e4.10}.

\noindent Pick now any small $\e>0$ and consider the function
$$
\underline u(x,y)=(1-\e)\biggl(1-e^{\lambda(x-\varphi_{M,N}(y))}\biggl)^+.
 $$
 We claim that $\underline u\leq u$ in ${\Gamma}_{M,N}\cap\{y\leq N\}$. Indeed, recall that we have
$$
 \lim_{x-x_{M,N}\to-\infty}u(x,N)=1,
$$
 this is a simple consequence of the uniform convergence of $u(x,.)$ to 1 as $x\to-\infty$ in the original variables, the fact that $(x_{M,N},N)$ is at the left of the leftmost point of $\Gamma$, and \eqref{e4.11}, and the fact that $\underline u\leq 1-\e$.
 Inside ${\Gamma}_{M,N}\cap\{y\leq N\}$, we have
 $$
 \begin{array}{rll}
& -\Delta u+c\partial_xu\\
=&\biggl(-\varphi_{M,N}''+\lambda{\varphi_{M,N}'}^2-\lambda-c\biggl)\lambda e^{\lambda(x-\varphi_{M,N}(y)}\\
 \leq&\biggl(\di\frac1{\sqrt{\inf(M,N)}}-4\lambda\mathrm{cotan}^2\beta+\lambda-c\biggl)\lambda e^{\lambda(x-\varphi_{M,N}(y)}\\
 \leq&0\ \hbox{as soon as $\lambda$ is small enough and $M,N$ large enough.}
 \end{array}
 $$
 Therefore, $\underline u\leq u$, so that we have, if we still call $\bar X$ the original point $\bar X$ translated by $\bar X_{M,N}$:
 $$
 u(\bar X)\geq 1-e^{-\lambda (M-\sqrt{\inf(M,N)})},
 $$
 which is the sought for estimate. \hfill$\Box$
 
 \smallskip
 \noindent{\bf Proof of Theorem \ref{t1.3}.} Proposition \ref{p4.1} proves that any limit of a sequence of translations $(u_n)_n$ in \eqref{e4.9} satisfies the assumptions of Theorem \ref{t4.1}. Therefore it is either a conical front $u_\alpha$ (Case 1), a tilted one-dimensional wave $\phi_0(y\cos\alpha-x\sin\alpha)$ (Case 2), or a tilted wave $\phi_0(y\cos\alpha+x\sin\alpha)$ (Case 3). We wish to prove that only Case 3 survives. Let us consider the set $T$ of turning points of $\Gamma$, that is, the set of all points 
 $(\varphi(y),y)$ such that $\varphi'(y)=0$. From the analyticity of $\Gamma$, $T$ is discrete: 
 $$
 T=\{(\varphi(y_n),y_n),\ n\in\NN\}.
 $$
 If we manage to prove that it is finite we are done because this excludes Case 1 trivially, and Case 2 because $u(x,y)\to 1$ as $x\to-\infty$, uniformly in $y$. So, assume that $T$ is infinite. We claim the existence of two  sequences of $T$, decreasing to $-\infty$, $(y_n^i)_n$, $i\in\{1,2\}$, such that 
 \begin{itemize}
 \item we have $y_{n+1}^1\leq y_n^2\leq y^1_n$ for all $n$,
 \item we have $\di \lim_{n\to+\infty}(y_n^2-y_{n+1}^1)=\di \lim_{n\to+\infty}(y_n^1-y_{n}^2)=+\infty$,
 \end{itemize}
  such that, if we set $x_n^i=\varphi(y_n^i)$ we have
$$
\lim_{n\to+\infty}u(x_n^i+x,y_n^i,y)=u_\alpha(x,y),
$$
uniformly on compact sets in $(x,y)$.
Indeed, there is $N_0\in \NN$ such that $T$ may be organised in disjoint clusters
 $$ 
 T=\bigcup_{n\in\NN}\Lambda_n,\ \hbox{$y\leq z$ if $(\varphi(y),y)\in\Lambda_{n+1}$, $(\varphi(z),z)\in\Lambda_n$},
 $$
 each $\Lambda_n$ having at most $N_0$ elements, and 
 $$
 \lim_{n\to+\infty}\biggl(\min_{(y,\varphi(y))\in\Lambda_n}y-\max_{(z,\varphi(z))\in\Lambda_n}z\biggl)=+\infty.
 $$
 This is because of the fact that $u_\alpha$ has a finite number of turning points and that the convergence of $u$ to its limits implies the convergence of the free boundaries in $C^1$ norms. We claim that $(x_n^1,y_n^1)$ and $(x_n^2,y_n^2)$ cannot be in two consecutive clusters, because of the orientation of $\mathcal{C}_\alpha$. So, if $\Lambda_n^i$ is the cluster of $(x_n^i,y_n^i)$, let $\Lambda^3_n$ be a cluster in-between. Let  $(x_n^3,y_n^3)$ be the leftmost point of $\Lambda_n^3$. From the definition of $\Lambda_n^3$, there is $R_n\to+\infty$ such that $(x_n^3,y_n^3)$ is the leftmost point of $\Gamma$ in $B_{R_n}(x_n^3,y_n^3)$. Let $u_\infty(x,y)$ be a limit of translations of $u$ with the sequence of points $(x_n^3,y_n^3)$, it does not converge to any of the functions prescribed by Theorem \ref{t4.1}, which is a contradiction. This finishes the proof of Theorem \ref{t1.3}. \hfill$\Box$
 \section{The model with two species}
 The first thing one must understand is \eqref{e1.1} in the truncated cylinder $\Sigma^L$, namely
 \begin{equation}
\label{e5.1}
\left\{
\begin{array}{rll}
-d\Delta v+c\partial_x v=&0\quad(x,y)\in\Sigma^L\cap\{v>0\}\\
\vert\nabla v\vert=&1\quad((x,y)\in\Gamma:=\partial\{v>0\}\cap\Sigma^L\\
\ \\
-Du''+cu'+1/\mu v_y=&0\quad\hbox{for $x\in\RR$, $y=0$}\\
dv_y(x,0)=&\mu u(x)-v(x,0)\quad\hbox{for $x\in\RR$}\\
\  \\
v(-\infty,y)=&1,\quad v(+\infty,y)=(1-y-L)^+\\
v(x,-L)=&1.
\end{array}
\right.
\end{equation}

\subsection{Exponential solutions}
\noindent System \eqref{e5.1}, linearised around 0, reads
\begin{equation}
\label{e5.2}
\left\{
\begin{array}{rll}
-d\Delta v+c\partial_x v=&0\quad(x,y)\in\Sigma^L\\
\ \\
-Du''+cu'+1/\mu v_y=&0\quad\hbox{for $x\in\RR$, $y=0$}\\
dv_y(x,0)=&\mu u(x)-v(x,0)\quad\hbox{for $x\in\RR$}\\
v(x,-L)=&0\\
v(-\infty,y)=&0.
\end{array}
\right.
\end{equation}
The solutions that decay to 0 as $x\to-\infty$ are looked for under the form 
$$(\phi_D^L(x),\psi_D^L(x,y))=e^{\alpha^L_D x}(1,\gamma\mathrm{sh}(\beta^L_D(y+L)),
$$
so that the exponents $\alpha^L_D$ and $\beta^L_D$ satisfy
\begin{equation}
\label{e5.3}
\left\{
\begin{array}{rll}
-d(\alpha^2+\beta^2)+c\alpha=&0\\
-D\alpha^2+c\alpha+\di\frac{d\beta\mu\mathrm{ch}(\beta L)}{\mathrm{sh}+d\beta\mathrm{ch}(\beta L)}=&0.
\end{array}
\right.
\end{equation} 
Once again we consider types of limits.

\noindent{\bf Case 1.} The limit $D\gg L\gg1$, $c$ bounded. We expect $\beta_D^L$ to go to 0 as $D\to+\infty$, so that 
$$
\frac{d\beta\mathrm{ch}(\beta L)}{\mathrm{sh}(\beta L)+d\beta\mathrm{ch}(\beta L)}\sim\frac1L.
$$
Then \eqref{e5.3} yields 
estimates of the form 
\begin{equation}
\label{e5.4}
\alpha^L_D\sim\sqrt{\frac\mu{ LD}},\quad \beta^L_D\sim\sqrt{\frac{c\sqrt\mu }{d\sqrt{DL}}}.
\end{equation}
And, once again, the estimate  may be pushed up to $c=o(\sqrt D)$.  

\noindent{\bf Case 2.}  The limit $D\gg L\gg1$, $c\gg\sqrt D$. This time we  expect $\beta_D^L$ to go to infinity as $D\to+\infty$, so that 
$$
\frac{d\beta\mathrm{ch}(\beta L)}{\mathrm{sh}(\beta L)+d\beta\mathrm{ch}(\beta L)}\to 1.
$$
We have
\begin{equation}
\label{e5.5}
\alpha_D^L\sim\biggl(\frac{c}{dD^2}\biggl)^{1/3},\quad \beta_D^L\sim\frac1\mu\biggl(\frac{c}{\sqrt D}\biggl)^{2/3}.
\end{equation}

\noindent{\bf Case 3.}  The limit $L\gg D\gg1$, $c\gg1$. We expect that $L\beta_D^L$ will go to infinity, so that $\di\frac{d\beta\mu\mathrm{ch}(\beta L)}{\mathrm{sh}(\beta L)+d\beta\mathrm{ch}(\beta L)}\to 1.$ In this setting, we have the estimates \eqref{e5.5}.

\subsection{Estimates on the velocity}
Let $(c_D^L,\Gamma_D^L,v_D^L)$ a solution to \eqref{e5.1}, we may infer its existence from a once again slight modification of Theorem 1.1 of \cite{CafR}.  We assume that it meets the line $\{y=0\}$ at the point $(0,0)$. 
\begin{thm}
\label{t5.1}
There is a constant $K>0$, independent of $L$ and $D$, such that
\begin{equation}
\label{e5.6}
c_D^L\leq K\sqrt D.
\end{equation}
Moreover there is $L_0>0$ such that we have
\begin{equation}
\label{e5.7}
\lim_{D\to+\infty}c_D^L=+\infty,
\end{equation}
uniformly in $L\geq L_0$.
\end{thm}
\noindent {\bf Proof.} That of \eqref{e5.6} is similar to that of Theorems \ref{t2.10}, one compares $(u_D^L,v_D^L)$ to 
$$
(\underline u_D^L(x),\underline v_D^L(x,y))=(\frac1\mu-\phi_D^L(x),1-\psi_D^L(x,y)).
$$
As for \eqref{e5.7}, the only point to prove is that, under the assumption that $(c_D^{L_0})_D$ is bounded, the width $L_0$ being large but fixed, then the leftmost point of $\Gamma_D^L$, denoted by $(\bar x_D^{L_0},\bar y_D^{L_0})$ satisfies
\begin{equation}
\label{e5.8}
\bar x_D^{L_0}\leq-\delta_0\sqrt D.
\end{equation}
Once this is proved, the proof of the theorem proceeds much as that of Theorem \ref{t2.1}, using in particular estimates \eqref{e5.4} for the linear exponentials. So, assume \eqref{e5.8} to be false, that. is, there are two diverging positive  sequences $(D_n)_n$
and $(x_n)_n$, and a positive constant $d_0$ such that
\begin{equation}
\label{e5.9}
\ \lim_{n\to+\infty}\frac{x_n}{\sqrt {D_n}}=0,\quad \min_{x\leq -x_n}\{y<0:\ (x,y)\in\Gamma_{D_n}^{L_0}\}\leq-d_0.
\end{equation}
By nondegeneracy, there is $\gamma>0$ such that 
\begin{equation}
\label{e5.9}
v(x,-\frac{d_0}2)\geq \gamma\ \hbox{for $x\leq x_n$}.
\end{equation}
Choose another sequence $(x_n')_n$ such that 
$$
\lim_{n\to+\infty}\frac{x_n}{x_n'}=\lim_{n\to+\infty}\frac{x_n'}{D_n}=0.
$$
Then we have
$$
\lim_{n\to+\infty}\Vert u_{D_n}^{L_0}\Vert_{L^\infty(-x_n',-x_n)}=0,
$$
just by integrating the ODE for $u_{D_n}^{L_0}$. But then, the Robin condition $\partial_y v_{D_n}^{L_0}+v_{D_n}^{L_0}=o(1)$ on $(-x_n',-x_n)\times\{0\}$, together with \eqref{e5.9} and the Hopf Lemma, yields the existence of $\gamma'>0$ such that 
$$
v_{D_n}^{L_0}(x,0)\geq \gamma'\ \hbox{for $-x_n'\leq x\leq x_n$.}
$$
The equation for $u_{D_n}^{L_0}$ becomes
$$
-(u_{D_n}^{L_0})''+c_{D_n}{L_0}(u_{D_n}^{L_0})'\leq\mu\un_{(-x_n,0)}(x)-\frac{\gamma'}2\un{(-x_n',-x_n)},
$$
as soon as $n$ is large enough. This implies
$
u_{D_n}{L_0}(-x_n')<0,
$ contradiction. \hfill$\Box$

\medskip
From then on, the rest of the study of the two species model parallels exactly that of the one species model. 
\section{Exponential convergence}
\noindent In this section, we assume that all the requirements on the coeffiients are fulfilled, and we drop the index ${}_D$ for the velocity $c$, the free boundary $\Gamma$ and the solution $u$. Theorem \ref{t1.4} is proved by deriving a differential inequlity for $\varphi$, exploiting  Theorem \ref{t1.2} and the fact that $\varphi'$ as a limit at infinity. This allows indeed to write the free boundary problem for $u$ in a suitable perturbative form, and translate the double Dirichlet and Neumann boundary condition into the sought for differential inequality.

\noindent In the whole section, the considerations will be rigorously identical for the one species model or the two species model, except the global estimate in Proposition \ref{p5.2} below, where the computations are slightly different - but left to the reader. Thus we will concentrate on the one species model (\ref{e1.2}).

\noindent Recall from \cite{ACF} that $\Gamma$ is an analytic graph. Theorem \ref{t1.2} readily implies
\begin{prop}
\label{p5.1}
We have 
$$\lim_{y\to-\infty}\varphi''(y)=\lim_{y\to-\infty}\varphi'''(y)=0.
$$
\end{prop}
This entails the following improvement of Proposition \ref{p4.2}.
\begin{prop}
\label{p5.2}
There is  $\rho>0$ such that , if $u(x,y)>0$ we have
$$
u(x,y)\geq 1-e^{-\rho\mathrm{dist}((x,y),\Gamma)}.
$$
\end{prop}
\noindent {\bf Proof.} For $\e>0$ consider $\psi_\e(y)$ smooth whose derivative $\psi_\e'$ satisfies
$$
\left\{
\begin{array}{rll}
\psi_\e(y)=\varphi'(y)\ &\hbox{if $y\leq-\di\frac1{\e^4}$}\\
\psi_\e'(y)=-\e\ &\hbox{if $-\di\frac1{\e^2}y\leq y\leq0$}\\
-\e\leq\psi_\e'(y)\leq0\ &\hbox{everywhere},\\
\vert \psi_\e''\vert\leq\e \ &\hbox{everywhere}.
\end{array}
\right.
$$
Note that this function exists due to Theorem \ref{t1.2} and Proposition \ref{p5.1}. For $\rho>0$ consider
$$
\bar u(x,y)=e^{\rho(x-\psi_\e(y))},
$$
it can be estimated by $e^{-\delta_\e\rho\mathrm{dist}((x,y),\Gamma)}$ for a suitable $\delta_\e>0$, because of Theorem \ref{t1.2} again. In the lower half plane we have
$$
-d\Delta \bar u+c\partial_x\bar u\\
=\biggl(-d(\rho-\psi_\e''+\rho(\psi_\e')^2)+c \biggl)\rho e^{\rho(x-\psi_\e(y))}
\geq0\ \hbox{if $\rho\in(0,c)$ and if $\e$ small enough.}
$$
On the line we have, on the same pattern:
$$
\biggl(-D\partial_{xx}+c\partial_x+\di\frac1\mu\partial_y\biggl)\bar u
=\biggl(-\rho D+c-\di\frac\e\mu \biggl)\rho e^{\rho x}
\geq0\ \hbox{if $\rho\in(0,\di\frac{c}D)$ and if $\e$ small enough.}
$$
And so, as soon as $\rho\in(0,\di\frac{c}D)$ and $\e>0$ is small enough, the function $\underline u:=(1-\bar u)^+$ is a  sub-solution to the equations for $u$ in the region $\{u>0\}$, moreover $\partial\{\bar u>0\}$ coincides $\Gamma$ sufficiently far in the lower half plane. The maximum principle implies the proposition. \hfill$\Box$    
\noindent{\bf Proof of Theorem \ref{t1.4}.} From now on, translate the origin so that we are in the following situation: $\Gamma\cap\RR^2_-$ is the graph
$$
\{(x,\phi(x)),\ x>0\},
$$
with $\phi'(x)<0$ and $\di\lim_{x\to+\infty}\phi'(x)=-\mathrm{tan}\alpha.$ Actually, the translation may be adjusted so that $\phi'(x)$ is as colse as we wish to $-\mathrm{tan}\alpha$, this will be quantified later. The following chain of transformations is then made.
\begin{enumerate}
\item Rotate the coordinates by the angle $\alpha$, so as to obtain the new set $(x,Y)$ given by
$$
X=x\cos\alpha-y\sin\alpha,\ \ Y=x\sin\alpha+y\cos\alpha.
$$
In this new system the free boundary may be written as
$Y=\psi(X),$ with, by Proposition \ref{p5.1}:
$$\lim_{X\to+\infty}\psi'(X)= \lim_{X\to+\infty}\psi''(X)= \lim_{X\to+\infty}\psi'''(X)=0.
$$
\item Straighten the free boundary by setting
$$
X'=X,\ Y'=\psi(X).
$$
In this new coordinate system, $u$ solves the over-determined problem
\begin{equation}
\label{e5.10}
\begin{array}{rll}
-\partial_{X'X'}u-(1+{\psi'}^2)\partial_{Y'Y'}u+\psi''\partial_{Y'}u-2\psi'\partial_{X'Y'}u&\\
+c\cos\alpha\partial_{X'}
u+(c_0-c\psi'\cos\alpha)\partial_{Y'}u=&0\quad{(X'>0,Y'<0)}\\
u(X',0)=0,\ \
\partial_{Y'}u(X',0)=&\di\frac1{\sqrt{1+{\psi'}^2(X')}}
\end{array}
\end{equation}
\item A standard compactness/uniqueness argument shows that
$$\lim_{X'\to+\infty}u(X',Y')=\phi_0(Y'),
$$
uniformly in $Y'<0$. As we will not need any additional change of coordinates, let us, for notational simplicity, revert to the initial notation
$$
x:=X',\ y:=Y'.
$$
The function $u(x,y)$ is this looked for under a perturbation of $\phi_0(y)$: $u(x,y)=\phi_0(y)+\tilde u(x,y)$. The free boundary condition writes
$$
\partial_yu(x,0)=-\frac{{\psi'}^2(x)}{\sqrt{1+{\psi'}^2(x)}}.
$$
The full PDE for $\tilde u$ will not be written, as we need another change of unknowns.
\item In order to transform \eqref{e5.10} into an over-determined problem with fixed boundary, we look for $\tilde u $ under the form
$$
\tilde u(x,y)=v(x,y)-\frac{{\psi'}^2(x)}{1+{\psi'}^2(x)+\sqrt{1+{\psi'}^2(x)}}\gamma_0(y),
$$
with $\gamma_0$ smooth, compactly supported, $\gamma(0)=0$, $\gamma'(0)=1$. As $\psi$ appears in the equations only in the form of $\psi'$, we set $h(x)=\psi'(x)$.
\end{enumerate}
The system for $(h,v)$ is thus
\begin{equation}
\label{e5.11}
\begin{array}{rll}
-\partial_{xx}v-(1+h^2)\partial_{yy}v+h'\partial_{y}v-&2h\partial_{xy}v+c\cos\alpha\partial_{x}v
+(c_0-ch\cos\alpha)\partial_{y}v\\
=&(h'+hc\cos\alpha)\phi_0'+{\mathcal R}[h](x,y)\quad(x>0,y<0)\\
v(x,0)=&0,\ 
\partial_{y}v(x,0)=0,
\end{array}
\end{equation}
the function ${\mathcal R}[h](x,y)$ being:
$$
\begin{array}{rll}
{\mathcal R}[h](x,y)=&\di{\biggl(\frac{h^2}{1+h^2+\sqrt{1+h^2}}\biggl)'' +\frac{h^2\gamma_0''(y)\sqrt{1+h^2}}{1+\sqrt{1+h^2}}-\frac{h^2h'\gamma_0'(y)}{1+h^2+\sqrt{1+h^2}}}\\
&\di{+2h\biggl(\frac{h^2}{1+h^2+\sqrt{1+h^2}}\biggl)'\gamma_0'(y)-c\cos\alpha\biggl(\frac{h^2}{1+h^2+\sqrt{1+h^2}}\biggl)'\gamma_0(y)}\\
&\di{-(c_0-hc\cos\alpha)\frac{h^2\gamma_0'(y)}{1+h^2+\sqrt{1+h^2}}+h^2\phi_0''(y)}
\end{array}
$$
While the expression of ${\mathcal R}[h](x,y)$ is utterly unpleasant, its structure is quite simple: it is quadratic in $h$ and its derivatives up to order 2, which are known to vanish at infinity. Finally, set
$$
{\mathcal S}[h,v](x,y)=h^2v_{yy}-h'v_y+2hv_{xy}+hc\cos\alpha v_y.
$$
By elliptic regularity, all derivatives of $v$ go to 0 as $x\to+\infty$, uniformly in $y$. Problem \eqref{e5.11} now reads
\begin{equation}
\label{e5.11}
\begin{array}{rll}
-\Delta v+c\cos\alpha v_x+c_0v_y
=&(h'+hc\cos\alpha)\phi_0'\\
&+{\mathcal R}[h](x,y)+{\mathcal S}[h,v](x,y)\quad(x>0,y<0)\\
v(x,0)=&0,\ 
\partial_{y}v(x,0)=0,
\end{array}
\end{equation}
We will use the estimate
\begin{equation}
\label{e5.15}
\vert {\mathcal R}[h](x,y)\vert+\vert {\mathcal S}[h,v](x,y)\vert\leq\theta(x,y)(\vert h'(x)\vert+\vert h(x)\vert)\un_{(0,1)}(y),
\end{equation}
 with $\theta(x,y)$  a positive function that tends to 0, uniformly in $y$, as $x\to+\infty$. It only remains to apply Proposition A1 of the appendix together with \eqref{e5.15}), with the data
 $$
 f(x)=h'(x)-h(x)c\cos\alpha,\ \ g(x,y)={\mathcal R}[h](x,y)+{\mathcal S}[h,v](x,y).
 $$
 From Proposition \ref{p5.2}, estimate (A1) holds, so that the function $h$ solves the differential inequality
 $$
 \vert h'(x)-h(x)c\cos\alpha\vert\leq\tilde\theta(x)(\vert h'(x)\vert+\vert h(x)\vert),
 $$
 the function $\tilde \theta(x)$ going to 0 as $x\to+\infty$. This implies the exponential estimate. \hfill$\Box$

\bigskip
\noindent {\bf Acknowledgement.} L.A. Caffarelli is supported by NSF grant DMS-1160802. The research of J.-M. Roquejoffre  has received funding from the ERC under the European Union's Seventh Frame work Programme (FP/2007-2013) / ERC Grant Agreement 321186 - ReaDi. He also acknowledges a long term visit at UT Austin during the academic year 2018-19, made possible by a delegation within the CNRS for that year, and a J.T. Oden fellowship. 
\vfill\eject
\section*{Appendix: compatibility relations for over-determined equations with right handside}
\noindent the goal of this section is to prove that the right handside of the solution of.an elliptic equation in the lower half plane, both with Dirichlet and Neumann condition, has to satisfy an integral equation.

\noindent {\bf Proposition A1.} {\it 
Consider $f(x)$ and $g(x,y)$ two smooth, bounded functions defined on $\RR_+$ (resp. $\RR_+\times\RR_-$). Set $e_0(y)=e^{c_0y}$. Assume the existence of $\rho>0$ such that 
$$
\vert g(x,y)\vert=O(e^{\rho y}),\ y\leq0,\ \hbox{uniformly in $x\in\RR_+$}.
\leqno{(A1)} 
$$
Let $u(x,y)$ solve 
$$
\left\{
\begin{array}{rll}
-\Delta u+c\cos\alpha u_x+c_0u_y=&f(x)e_0(y)+k(x,y)\ \ (x>0,\ y\leq0)\\
u(x,0)=u_y(x,0)=&0,
\end{array}
\right.
\leqno{(A2)}
$$
in the classical sense. Also assume that $u$ and its derivatives satisfy the estimate (A1). There is a continuous function $K(x,y)$, defined on $[2,+\infty)\times\RR_-$, such that, for all $\e>0$ we have
$$
\vert K(x,y)\vert=O(e^{(\beta-\rho-\e)y-\rho x}),\ \ x\geq2,\ y\leq0,
$$
and such that 
$$
f(x)=\int_{-\infty}^0K(x,y)g(x,y)dy,\ \ x\geq2.
$$
}

\noindent{\bf Proof.} First, transform (A2) into a problem on the whole line by setting 
$$
v(x,y)=\gamma(x)u(x,y),
$$
where $\gamma$ is a smooth function, supported on $(1,+\infty)$, and equal to 1 on $[2,+\infty)$. the equation for $v$ is thus 
$$
\left\{
\begin{array}{rll}
-\Delta v+c\cos\alpha v_x+c_0v_y=&\gamma(x0f(x)e_0(y)+k(x,y)\ \ (x\in\RR,\ y\leq0)\\
v(x,0)=v_y(x,0)=&0,
\end{array}
\right.
\leqno{(A3)}
$$
with
$$
k(x,y)=\gamma(x)g(x,y)+\gamma''(x)u(x,y)+2\gamma'(x)u_x(x,y).
$$
In other words, the right handside of (A3) equals that of (A2) as soon as $x\geq2$.

\noindent A second observation is that the general Dirichlet problem 
$$
\left\{
\begin{array}{rll}
-\Delta u+c\cos\alpha u_x+c_0u_y=&F(x,y)\ \ (x>0,\ y\leq0)\\
u(x,0)=&0\\
u(0,y)=&u_0(y)
\end{array}
\right.
\leqno{(A3)}
$$
is well-posed as soon as the Dirichlet datum $u_0(y)$  and the right handside $F(x,y)$ satisfy the estimate  (A1), with any $\rho<c_0$. Indeed, the classical change of unknown $u(x,y)=e^{\rho y}v(x,y)$ changes (A3) into an elliptic equation the the zero order term
$\rho c_0-\rho^2$. Thus it is enough to assume tht $g$ and $k$ are compactly supported in $y$, in order to perform Fourier transforms. An easy density argument then allows to treat general data satisfying (A1).
 
\noindent Set $w(x,y)=e^{(c\cos\alpha x+c_0y)/2}v(x,y)$, the equation for $w$ is now
$$
\left\{
\begin{array}{rll}
-\Delta w+\beta^2w=&F(x)e_0(y)+G(x,y)\ \ (x\in\RR,\ y\leq0)\\
w(x,0)=w_y(x,0)=&0.
\end{array}
\right.
\leqno{(A4)}
$$
The new data are 
$$
F(x)=e^{-c\cos\alpha/ x/2}\gamma(x)f(x),\ G(x,y)=e^{-(c\cos\alpha x-c_0y)/2}k(x,y),\ E_0(y)=e^{c_0y/2},
$$
and we have set
$$\ \beta^2=\frac{c_0^2+c^2\cos^2\alpha}4.
$$
Thanks to both Dirichlet and Neumann conditions, the function $w(x,y)$ may be extended evenly over the whole plane $\RR^2$, where it solves the same equation as (A4), with $G$ and $E_0$ also extended evenly. The new unknown and data are still denoted by $w$, $E_0$ and $G$. Let $\hat w(\xi,\zeta)$ be the Fourier transform of $w$, as well as $\hat F(\xi)$, $\hat E_0(\zeta)$ and $\hat G(\xi,\zeta)$ the Fourier transforms of the data. Thus we have
$$
\hat w(\xi,\zeta)=\frac{\hat F(\xi)\hat E_0(\zeta)+\hat G(\xi,\zeta)}{\beta^2+\xi^2+\zeta^2}.
$$
Now, the Dirichlet condition for $w$ entails $\di\int_{\RR}\hat w(\xi,\zeta)d\zeta=0$, that is
$$
\hat F(\xi)=-\biggl(\int_{\RR}\frac{\hat E_0(\zeta)}{\beta^2+\xi^2+\zeta^2}d\zeta\biggl)^{-1}\int_{\RR}\frac{\hat G(\xi,\zeta)}{\beta^2+\xi^2+\zeta^2}d\zeta.
$$
From Plancherel's Theorem, we have 
$$
\int_{\RR}\frac{\hat E_0(\zeta)}{\beta^2+\xi^2+\zeta^2}d\zeta\propto\frac1{c_0/2+\sqrt{\beta^2+\xi^2}},
$$
and, if we still denote by $\hat G(\xi,y)$ the partial Fourier transform of $G$ with respect to $x$, we have, again by Plancherel's theorem
$$
\int_{\RR}\frac{\hat E_0(\zeta)}{\beta^2+\xi^2+\zeta^2}d\zeta\propto\int_{-\infty}^{+\infty}\hat G(\xi,y)e^{-\vert y\vert\sqrt{\beta^2+\xi^2}}dy.
$$
Inverting the Fourier transform in $x$ yields, by the residue theorem:
$$
F(x)=\propto e^{-c\cos\alpha\vert x\vert/2}\int_{\RR^2}G(x,y)e^{ix\xi-\vert y\vert\sqrt{\beta^2+\xi^2}}d\xi dy.
$$
Shifting the integration line in $\xi$ from $\RR$ to $\RR+i\rho$ yields, for $x>0$:
$$
\int_{\RR}e^{ix\xi-\vert y\vert\sqrt{\beta^2+\xi^2}}d\xi=e^{-\rho x}\int_{\RR}e^{ix\xi-\vert y\vert\sqrt{\beta^2+\xi^2-\rho^2+2i\rho\xi}}d\xi
$$
and $\sqrt{\beta^2+\xi^2-\rho^2+2i\rho\xi}$ has nonzero real part. This entails the proposition. \hfill$\Box$
\noindent 
{\footnotesize
 
}

\end{document}